\documentclass[12pt]{article}
\usepackage{mathrsfs}
\usepackage{tipa}
\usepackage{amssymb}
\usepackage{amsfonts}

\usepackage{graphicx}
\usepackage{latexsym,amsmath,amssymb,amsfonts,amsthm}

\newtheorem{theorem}{Theorem}[section]
\newtheorem{lemma}[theorem]{Lemma}

\newtheorem{corollary}[theorem]{Corollary}

\begin{document}
\setcounter{page}{1}
\title{The dual foliation of some singular Riemannian foliations}
\author{Yi Shi}
\date{}
\protect
\maketitle ~~~\\[-5mm]
{ \footnotesize  School of Mathematical Sciences, Shanxi University, Taiyuan 030006, China,
e-mail: shiyi@sxu.edu.cn}\\[2mm]
{\bf Abstract:} \noindent In this paper, we use the methods of
subriemannian geometry to study the dual foliation of the singular
Riemannian foliation induced by isometric Lie group actions on a
complete Riemannian manifold $M$. We show that under some conditions, the dual foliation
has only one leaf.\\[2mm]
{\bf MSC(2010) Subject Classification:} 53C12, 53C17.

\markright{\sl\hfill The dual foliation of some singular Riemannian foliations\hfill}

\section{Introduction}

\renewcommand{\thesection}{\arabic{section}}
\renewcommand{\theequation}{\thesection.\arabic{equation}}
\setcounter{equation}{0}

We recall some basic notions about singular Riemannian foliations,
for further details we refer the readers to \cite{abt,wa,mol,wi}. A \emph{singular
Riemannian foliation} $\mathcal{F}$ on a Riemannian manifold $M$ is a decomposition of $M$
into smooth injectively immersed submanifolds $L(p)$, called leaves,
such that it is a singular foliation and  any geodesic starting
orthogonally to a leaf remains orthogonal to all leaves it
intersects. Such a geodesic is called a \emph{horizontal geodesic}.
For all $p \in M$, we denote by $H_{p}$ the orthogonal complement to
the tangent space $T_{p}(L(p))$, and call it the horizontal space at
$p$. If all the leaves have the same dimension, then $\mathcal{F}$ is called a \emph{regular
Riemannian foliation}.

A curve is called \emph{horizontal} if it meets the leaves of
$\mathcal {F}$ perpendicularly. In \cite{wi}, Wilking associates to
a given singular Riemannian foliation $\mathcal{F}$ the so-called
 \emph{dual foliation} $\mathcal{F}^{\#}$. The \emph{dual
leaf} through a point $p \in M$ is defined as all points $q \in M$
such that there is a piecewise smooth, horizontal curve from $p$ to
$q$. We denote by $L^{\#}_{ p}$ the dual leaf through $p$. Wilking proved
that when $M$ has positive curvature, the horizontal connectivity
holds on $M$, i.e. the dual foliation has only one leaf. Wilking also used the theory of
dual foliations to show that the Sharafutdinov projection is
smooth. For more applications of dual foliations, the reader is referred to \cite{fgt,gz,gw1,gw2,ly,wi}.

Let $(M, \mathcal{F})$ be a singular Riemannian foliation $\mathcal{F}$ on Riemannian manifold $M$. Applying Wilking's theorem, we can get the horizontal connectivity by assuming that $M$ has positive curvature. However, in this paper we are interested in getting the horizontal connectivity by applying the methods of subriemannian geometry. For a more extensive exposition of subriemannian geometry, we refer the readers to \cite{br,cc,mon}.

Recall that a \emph{subriemannian geometry} on a manifold $M$ consists of a distribution, which is to say a vector subbundle $H\subset TM$ of the tangent bundle of $M$, together with a fiber inner-product on $H$. We call $H$ the horizontal distribution. A curve on $M$ is called horizontal if it is tangent to $H$.

Let $H$ be the horizontal distribution on $M$, the Lie brackets of
vector fields in $H$ generate the flag
$$H=H^{1}\subset H^{2}\subset \cdots \subset H^{r}\subset \cdots\subset TM$$
with$$H^{r+1}=H^{r}+[H,H^{r}] \ \  for \ \  \textit{r}\geq 1$$
where$$[H,H^{r}]=span\{[X,Y]:X\in H,Y\in H^{r}\}$$

At a point $p\in M$, this flag  gives a flag of subspaces of
$T_{p}M$:
$$H_{p}=H^{1}_{p}\subset H^{2}_{p}\subset \cdots \subset H^{r}_{p}\subset \cdots\subset T_{p}M$$

We say that $H$ is
\emph{bracket generating} at p if there is an $r\in \mathbb{Z}^{+}$
such that $H^{r}_{p}= T_{p}M$, and  $H$ is \emph{bracket generating}
if for all $x\in M$ there is an $r(x)\in \mathbb{Z}^{+}$ such that
$H^{r(x)}_{ x }= T_{x}M$. The smallest integer $r$ such that $H^{r}_{p}=T_{p}M$ is called
the \emph{step} of the distribution at $p$. \\[2mm]
\noindent $\mathbf{ Chow's\ \ Theorem}$ (\cite{ch}). Let M be a
connected manifold and $H \in TM$ be a bracket generating
distribution, then the set of points that can be connected to $p\in
M$ by a horizontal
path coincides with $M$.\\[2mm]
\indent When Chow's condition fails on some subset of $M$, sometimes
the horizontal connectivity also fails (\cite{cc}, P. 82).

Suppose a compact connected Lie group $G$ acts isometrically on manifold $M$. Throughout the paper, every action will be assumed to be \emph{effective}. It is well known that the action of $G$ induces a singular Riemannian foliation $(M, \mathcal{F})$ on $M$. Let $M'$ denote the union of all principal orbits in $M$. An basic fact is that the restricted foliation $(M', \mathcal{F}|_{M'})$ is a regular Riemannian foliation on $M'$. Denote by $H'$ the collection of horizontal spaces $H_{x}$ with $x\in M'$, then $H'$ is the horizontal distribution of $\mathcal{F}|_{M'}$ on $M'$. Contrast to Chow's Theorem, to connect any two points on
$M$ by a horizontal curve, we only need to assume that $H'$ is bracket-generating at one point:

\begin{lemma}\label{lem1}
Let $G$ be a compact connected Lie group, acting isometrically on a complete Riemannian manifold $M$. If $H'$ is bracket-generating at some point $p\in M'$, then the dual foliation has only one leaf. If $M$ is compact, then there exists a constant $C=C(M, G)$ such that any two points of $M$ can be connected by a horizontal curve of length $\leq C$.
\end{lemma}

Let $G \times M \rightarrow M$ denote an isometric action by principal orbits of a compact Lie
group $G$ on a complete Riemannian manifold $M$. The orbit space $B:= M/G$ is also a manifold, and inherits a quotient metric from
$M$ for which $\pi: M \rightarrow B$ becomes a Riemannian submersion. $\pi$ is then said to be a \emph{homogeneous submersion} or \emph{fibration}.

The following corollary of Lemma 1.2 gives a generalization of Theorem 3.1 in \cite{gw1}.

\begin{corollary}\label{cor}
Let $\pi: M \rightarrow B= M/G$ be a homogeneous $G$-fibration with compact connected Lie group $G$. If the
horizontal distribution $H$ is bracket-generating at some point
$p\in M$, then $\pi: M \rightarrow B$ is a principal $G$-bundle.
\end{corollary}

Using Lemma 1.2 we can also get the following:

\begin{theorem}\label{T1}
Suppose circle $T^{1}$ acts isometrically on a complete manifold $M$. If $M$ is simply connected, then the dual foliation has only one leaf.
If $M$ is compact, then there exists a constant $C=C(M, T^{1})$ such that any two points of $M$ can be connected by a horizontal curve of length $\leq C$.
\end{theorem}

Notice that Wilking's horizontal connectivity theorem doesn't hold for nonnegatively
curved manifolds in general, but as showed in \cite{agw}, the horizontal connectivity phenomenon often occurs for nonnegatively curved spaces. The main result of this paper continues this investigation:

\begin{theorem}\label{Tk} Suppose torus $T^{k}$ acts isometrically on a complete nonnegatively curved
Riemannian manifold $M$. Then either the dual foliation has only one
leaf, or else $M$ locally splits. \end{theorem}

As an application of Theorem 1.3 and Theorem 1.4, the following theorem gives a generalization of Theorem 4 in \cite{agw}.

\begin{theorem}
Suppose torus $T^{k}$ acts isometrically on a complete, simply
connected Riemannian manifold $M$. Then the dual foliation has only one
leaf if one of the following holds:

(a) the action of $T^{k}$ is free.

(b) $M$ has nonnegative curvature.
\end{theorem}

\section{Bracket-generating at one point}
\renewcommand{\thesection}{\arabic{section}}
\renewcommand{\theequation}{\thesection.\arabic{equation}}
\setcounter{equation}{0}
\setcounter{theorem}{0}
In this section, we use the methods of subriemannian geometry to give a proof of Lemma 1.1. The proof of Corollary \ref{cor} follows as an application of Lemma 1.1.\\[2mm]
\noindent {\it Proof of Lemma \ref{lem1}} Suppose that $H'$ is bracket-generating at some point $p\in M'$. Since $M'$ is open and dense in $M$, by the continuity of Lie brackets, $H'$ is also bracket-generating in a neighbourhood $U_{p}\subset M'$ of $p$. By the Ball-box theorem (\cite{mon}, P. 29) we know that there is a positive constant
$\epsilon_{0}$ and a box neighborhood $\mathrm{Box}(\epsilon_{0}, p)\subset U_{p}$ of $p$, such that any
point in $\mathrm{Box}(\epsilon_{0}, p)$ can be connected to $p$ by a horizontal curve of length $\leq \epsilon_{0}$.

We claim that $G(p)\subset L^{\#}_{p}$, where $G(p)$ is the orbit of
$p$. For any $q\in G(p)$, assume $q=gp$ for some $g\in G$. Denote the
differential of $g$ by $g_{\ast}$, for any two vector fields $X$,
$Y\in H'$, we have $g_{\ast}[X,Y]=[g_{\ast}X,g_{\ast}Y]$. Assume that
$H'$ is bracket-generating at $p\in M'$
with step $r$, i.e., $(H')_{p}^{r}=T_{p}M$. Then by
$(H')_{q}^{r}=g_{\ast}((H')_{p}^{r})=g_{\ast}(T_{p}M)=T_{q}M$, $H'$ is also
bracket-generating at $q$ with step $r$. Since $g_{\ast}$ is a
isometric map, there is a same size box neighborhood
$\mathrm{Box}(\epsilon_{0}, q)$ of $q$, as that
of $p$, such that any point in $\mathrm{Box}(\epsilon_{0}, q)$ can
be connected to $q$ by a horizontal curve of length $\leq
\epsilon_{0}$. Since $G(p)$ is compact connected, we can cover
$G(p)$ by finitely many successive box neighborhoods. Denote the
number of these covering boxes by $N_{0}$. Then any two points $x, y\in G(p)$ can
be connected by a horizontal curve with length $\leq
2N_{0}\epsilon_{0}$. This proves that $G(p)\subset L_{p}^{\#}$.

For any $x \notin G(p)$, since $G(p)$ is compact, there is a minimal
geodesic $\gamma:[0,1]\rightarrow M$ from $x$ to $G(p)$. Thus $\gamma$ must
be a horizontal geodesic.  By $\gamma(1)\in G(p)\subset L_{p}^{\#}$
, we get $x\in L_{p}^{\#}$. This proves that $L_{p}^{\#}=M$, i.e.,
the dual foliation has only one leaf.

If $M$ is compact, then for any $x_{1},x_{2}\in M$, we have two minimal geodesics $\gamma_{i}:[0,1]
\rightarrow M$ from $\gamma_{i}(0)=x_{i}$ to $G(p)$ for $i=1,2$. Thus both $\gamma_{1}$ and $\gamma_{2}$ are horizontal geodesics. We have proved that $\gamma_{1}(1)$ and $\gamma_{2}(1)$ can be connected by a horizontal curve of length $\leq
2N_{0}\epsilon_{0}$. Denote by $diam(M)$ the diameter of $M$, since length$(\gamma_{i})\leq diam(M)$ for
$i=1,2$, $x_{1}$ and $x_{2}$ can be connected by a horizontal curve of length $\leq C(M,
G):= 2diam(M)+2N_{0}\epsilon_{0}$. $\hfill \square$\\

\noindent {\it Proof of Corollary \ref{cor}} For $p\in M$, let $G_{p}$ be the isotropy group of the action of $G$ at $p$. By Lemma \ref{lem1} we know that $L_{p}^{\#}=M$, then by \cite{wi} one can connect any two points of $M$ by a piecewise horizontal geodesic. In \cite{gw1} it was shown that, if $\gamma$ is a horizontal geodesic with respect to a homogeneous fibration of a manifold $M$, then the isotropy groups
of the action coincide for every point of $\gamma$. Hence, we get that the isotropy group of any point must coincide with $G_{p}$. Since the action of $G$ is effective, we have that $G_{p}$ is the identity, so the action of $G$ is free, and $\pi: M \rightarrow B$ is a principal $G$-bundle.  $\hfill \square$

\section{Torus action and horizontal connectivity}
\renewcommand{\thesection}{\arabic{section}}
\renewcommand{\theequation}{\thesection.\arabic{equation}}
\setcounter{equation}{0}
\setcounter{theorem}{0}

In this section, we study the dual foliation of the singular Riemannian foliation induced by isometric torus actions and prove Theorem 1.3, Theorem 1.4 and Theorem 1.5.\\[2mm]
\noindent {\it Proof of Theorem \ref{T1}} If there is a singular orbit
$x=T^{1}(x)$ on $M$, then the horizontal space at $x$ is
$H_{x}=T_{x}M$. Since the exponential map $\exp_{x}:
H_{x}\rightarrow M$ is onto, we get that $L_{x}^{\#}=M$. Thus we can
assume all the orbits are regular, and we get a codimension one
globally horizontal distribution $H$ on $M$. If $H$ is not
integrable at some point $p\in M$, then $H$ is bracket-generating at
$p$. By Lemma \ref{lem1} we get that $L_{p}^{\#}=M$. So it is enough to
show that there is no globally integrable horizontal distribution
$H$ on $M$.

If $H$ is integrable on $M$, since $M$ is simply connected, by Theorem 1.5 in \cite{at} there are no exceptional orbits on $M$, so $T^{1}$ acts
freely on $M$, it follows that the orbit space $B := M/T^{1}$ is
also a manifold, and inherits the quotient metric from $M$ for which
$\pi: M \rightarrow B$ becomes a Riemannian submersion with fiber
$F=T^{1}$. By the long exact homotopy sequence of the fibration
$\pi$, $B$ is simply connected since $M$ is. We claim that $M$ is
homeomorphic to $B\times T^{1}$.

Since $H$ is integrable, by the Frobenius theorem, for any $p\in M$,
$L_{p}^{\#}$ is the integral leaf of $H$ that contains $p$. Since
$\pi$ is a locally trivial bundle, $\pi\mid_{L_{p}^{\#}}:
L_{p}^{\#}\rightarrow B$ is a covering map and therefore a
diffeomorphsm, $B$ being simply connected. Hence $L_{p}^{\#}$ is a section of $\pi: M \rightarrow B$. Since it was proved already
that $\pi: M \rightarrow B$ is a principal bundle, the existence of a global section
implies immediately that $M$ is homeomorphic to $B\times T^{1}$. However, this is a contradiction
since $M$ is simply connected.

We have proved the horizontal connectivity on $M$,
assume now that $M$ is compact, denote by $diam(M)$ the original
diameter of $M$, by $diam_{H}(M)$ the horizontal diameter of $M$
defined by horizontal metric. If there is a singular orbit
$x=T^{1}(x)$ on $M$, then $L_{x}^{\#}=M$, thus $diam_{H}(M)\leq
2diam(M)$. If all the orbits are regular, then we get a codimension
one globally horizontal distribution $H$ on $M$. Since $H$ is
bracket-generating at some point $p\in M$,
by Lemma \ref{lem1}, we get that $diam_{H}(M)\leq C=C(M, T^{1})$. $\hfill \square$\\[2mm]
The following lemma will be used in the proof of Theorem \ref{Tk}:

\begin{lemma}\label{S}
Let $G$ be a compact connected Lie group, acting isometrically on a complete Riemannian manifold $M$. For any $p\in M$, define the subset $S\subset G$ via $S:=\{g\in G\ |\ gp\in L_{p}^{\#}\cap G(p)\}$. Then $S$ is a subgroup of $G$ and the action of $S$ on $M$ preserves $L_{p}^{\#}$.
\end{lemma}

\begin{proof}To show that $S$ is a subgroup of $G$, we need to prove that $g_{1}^{-1}g_{2}\in S$ for any $g_{1}, g_{2}\in S$. It is enough to show that $g_{1}^{-1}g_{2}(p)\in L_{p}^{\#}$. By \cite{wi} one can connect $p$ and $g_{i}(p)$ by a piecewise horizontal geodesic $\gamma_{i}$, $i=1,2$. It follows that one can connect $g_{1}^{-1}p$ and $p=g_{1}^{-1}(g_{1}p)$ by the piecewise horizontal geodesic $g_{1}^{-1}\circ\gamma_{1}$, thus $g_{1}^{-1}p\in L_{p}^{\#}$. Also, one can connect $g_{1}^{-1}p$ and $g_{1}^{-1}(g_{2}p)$ by the piecewise horizontal geodesic $g_{1}^{-1}\circ\gamma_{2}$.  Since $g_{1}^{-1}p\in L_{p}^{\#}$, $g_{1}^{-1}g_{2}(p)\in L_{p}^{\#}$ as desired.

We claim that the action of $S$ on $M$ preserves $L_{p}^{\#}$. For any $g\in S$ and $q\in L_{p}^{\#}$, we need to show that $gq\in L_{p}^{\#}$. Since one can connect $p$ and $q$ by a piecewise horizontal geodesic $\gamma$, one can connect $gp$ and $gq$ by the piecewise horizontal geodesic $g\circ\gamma$. Since $g\in S$, by definition $gp\in L_{p}^{\#}$. Thus $gq\in L_{p}^{\#}$ as claimed.
\end{proof}

\noindent {\it Proof of Theorem \ref{Tk}} Suppose the dual foliation has more than one leaf. Since $T^{k}$ acts transitively
on the space of dual leaves, all dual leaves have the same dimension. Then by theorem 2 and 3 in \cite{wi} these dual leaves form a regular Riemannian foliations $\mathcal{F}^{\#}$. Let $M'$ denote the union of all principal orbits in $M$, and let $p\in M'$. Since the torus is an abelian group, the isometric action of $T^{k}$ induces $k$ commuting Killing fields on $M$. Consider a Killing field $X$ which is perpendicular to $L_{p}^{\#}$ at $p$, we will prove that $X$ is perpendicular to $L_{p}^{\#}$ for all $x\in L_{p}^{\#}$. To prove this we may assume that $x\in L_{p}^{\#}\cap M'$ since $M'$ is open and dense in $M$. Denote the normal space of $L_{p}^{\#}$ at $x$ by $N_{x}(L_{p}^{\#})$. We first plan to show that $X(q) \in N_{q}(L_{p}^{\#})$ for all $q\in L_{p}^{\#}\cap T^{k}(p)$.

Denote the one-parameter subgroup of $T^{k}$ corresponding to $X$  by $\phi(s), s\in \mathbb{R}$. Define the subset $S\subset T^{k}$ via $S:=\{g\in T^{k}\ |\ gp\in L_{p}^{\#}\cap T^{k}(p)\}$, by Lemma \ref{S} the action of $S$ on $M$ preserves $L_{p}^{\#}$. Let $q\in L_{p}^{\#}\cap T^{k}(p)$, by definition $q=gp$ for some $g\in S$. Since $T^{k}$ is abelian, $g\circ \phi(s)(p)=\phi(s)(q)$. Denote the differential of $g$ by $g_{\ast}$, then $g_{\ast}(X(p))=X(q)$. Since $X(p)\in N_{p}(L_{p}^{\#})$ and the action of $S$ preserves $L_{p}^{\#}$, we get that $X(q)=g_{\ast}(X(p))\in N_{q}(L_{p}^{\#})$.

For any $x\in L_{p}^{\#}\cap M'$, we can find some $q\in L_{p}^{\#}\cap T^{k}(p)$ and a piecewise horizontal geodesic $\gamma$ in $M'$ connect $x$ and $q$. Since $X(q) \in N_{q}(L_{p}^{\#})$, by the proof of Corollary 8 in \cite{wi} $X$ is perpendicular to $L_{p}^{\#}$ along $\gamma$. Thus $X$ is perpendicular to $L_{p}^{\#}$ for all $x\in L_{p}^{\#}\cap M'$, which in turn shows $X$ is perpendicular to $L_{p}^{\#}$.

Denote the codimension of $L_{p}^{\#}$ by $n$, then $n\leq k$. We can choose $n$ Killing fields $X_{1},\cdots, X_{n}$ such that $\{X_{1}(p), \cdots, X_{n}(p)\}$ is an orthonormal basis of $N_{p}(L_{p}^{\#})$. We claim that the normal bundle of $L_{p}^{\#}$ is spanned by $\{X_{1},\cdots, X_{n}\}$.

For any $x\in L_{p}^{\#}$, we can find some $q\in L_{p}^{\#}\cap T^{k}(p)$ and an horizontal geodesic $c$ from $q=c(0)$ to $x=c(1)$. It is clear that $\{X_{1}(q), \cdots, X_{n}(q)\}$ is an orthonormal basis of $N_{q}(L_{p}^{\#})$. Set $W=span_{\mathbb{R}}\{X_{1},\ldots,X_{n}\}$ and $X(t)=X(c(t))$ for $X\in W$, $t\in [0,1]$, let $W(t)=\{X(t)|X\in W\}$. We will show that $X(t)\neq 0$ for any non-zero $X\in W$ and $t\in [0,1]$, thus $\dim(W(t))\equiv n$. Indeed, if $X(t_{0})=0$ for some $X:=\sum_{i\leq n}a_{i}X_{i}\in W$ and $t_{0}\in [0,1]$, then $X$ can be written also as a variation of $c$ by horizontal geodesics with a fixed value $c(t_{0})$ at time $t_{0}$.  Therefore $X$ is everywhere tangential to $L_{p}^{\#}$, so $X(0)=\sum_{i\leq n}a_{i}X_{i}(q)\perp N_{q}(L_{p}^{\#})$. But then $X(0)=0$ and $a_{i}=0$ for every $i\leq n$, hence $X=0$.

Since we have proved that any $X\in W$ is perpendicular to $L_{p}^{\#}$, now by $\dim(W(t))\equiv n$, we get that the normal bundle of $L_{p}^{\#}$ is spanned by $n$ linearly independent commuting Killing fields $X_{1},\cdots, X_{n}$. Thus the horizontal distribution with respect to $\mathcal{F}^{\#}$ is integrable, then by Theorem 1.3 in \cite{wal} $M$ locally splits. $\hfill \square$\\

\noindent {\it Proof of Theorem 1.5} (a) We first consider the case that the action of $T^{k}$ is free. We get a sequence of principal $T^{1}$-bundles
$$M=M_{0}\xrightarrow{\pi_{1}} M_{1}\xrightarrow{\pi_{2}}\cdots M_{k-1}\xrightarrow{\pi_{k}}
M_{k}=M/T^{k},$$ so $M_{i}$ is simply connected for any $0\leq i\leq k$. Notice that for any $0\leq
i\leq k-1$, there is a subgroup $T^{k-i}$ of $T^{k}$ acts freely on $M_{i}$, so this action induces a
Riemannian foliation $\mathcal {F}_{i}$ on $M_{i}$. For any $p_{i}\in M_{i}$,
denote by $L_{p_{i}}^{\#}$ the dual leaf of $\mathcal {F}_{i}$ that
through $p_{i}$. We claim that if $L_{p_{i}}^{\#}=M_{i}$ for some
$p_{i}\in M_{i}$ and $1\leq i\leq k-1$, then
$L_{p_{i-1}}^{\#}=M_{i-1}$ for any $p_{i-1}\in \pi^{-1}_{i}(p_{i})$.
This in fact shows $L_{p}^{\#}=M$ for any $p\in M$.

By Theorem \ref{T1} we get $L_{p_{k-1}}^{\#}=M_{k-1}$ for any
$p_{k-1}\in M_{k-1}$. Now suppose $L_{p_{i}}^{\#}=M_{i}$ for some $p_{i}\in M_{i}$ and $1\leq
i\leq k-1$. Assume to the contrary that $L_{p_{i-1}}^{\#}\neq M_{i-1}$
for some $p_{i-1}\in \pi^{-1}_{i}(p_{i})$. By \cite{wi}
$L_{p_{i-1}}^{\#}$ is a complete smooth immersed submanifold of
$M_{i-1}$ of dimension $\dim(L_{p_{i-1}}^{\#})\leq \dim(M_{i})$, it
is easy to see that $\pi_{i}\mid_{L_{p_{i-1}}^{\#}}:
L_{p_{i-1}}^{\#}\rightarrow L_{p_{i}}^{\#}=M_{i}$ is a smooth
surjective map, thus $\dim(L_{p_{i-1}}^{\#})= \dim(M_{i})$. Thus by the implicit function theorem
$\pi_{i}\mid_{L_{p_{i-1}}^{\#}}: L_{p_{i-1}}^{\#}\rightarrow
L_{p_{i}}^{\#}=M_{i}$ is a locally diffeomorphsm, and therefore a
diffeomorphsm, $M_{i}$ being simply connected. Now using the same
method as that of Theorem \ref{T1}, we get that $M_{i-1}$ is homeomorphic
to $M_{i}\times T^{1}$, which is a contradiction since $M_{i-1}$ is
simply connected.

(b) We next consider the case that $M$ has nonnegative curvature. Suppose the dual foliation has more than one leaf, then these dual leaves form a regular Riemannian foliations $\mathcal{F}^{\#}$. Denote by $H^{\#}$ the horizontal distribution with respect to $\mathcal{F}^{\#}$. Then by Theorem  \ref{Tk} $M$ locally splits. Since $M$ is simply connected, $M$ globally splits. Then $M$ is isometric to the Riemannian product $L_{p}^{\#}\times A$, where $L_{p}^{\#}$ is a dual leaf and $A$ is an integral manifold of $H^{\#}$. Therefore $L_{p}^{\#}$ is not dence. By Proposition 9.1 in \cite{wi} there is a closed subgroup $T^{n}\subsetneq T^{k}$ and a Riemannian submersion $\sigma:M\rightarrow T^{k-n}=T^{k}/T^{n}$, where the fibers of $\sigma$ are closures of leaves of $\mathcal{F}^{\#}$. By the long exact homotopy sequence of $\sigma$, $T^{k-n}$ is simply connected since $M$ is. This is a contradiction since $T^{k-n}$ is not simply connected. $\hfill \square$\\

\noindent$\mathbf{Acknowledgement.}$ The author thanks Martin Kerin for helpful comments on a previous version of this paper.

 \end{document}